\input ./style/arxiv-general.cfg
\documentclass[aos,MSNbibl,dvips]{arximspdf}
\makeatletter
   \@ifpackageloaded{graphicx}{}{\usepackage{graphicx}}
\makeatother
\usepackage{mathbh}
\usepackage{dcolumn}

\usepackage{graphicx}
\usepackage{epstopdf}

\doi{10.1214/15-AOS1381}% Updated by VTEXPTS2LaTeX.exe, 19.10.2015
\volume{44}
\issue{2}
\pubyear{2016}
\firstpage{660}
\lastpage{681}
\docsubty{FLA}

\makeatletter
\def\sfrac#1#2{#1/#2}
\def\vfrac#1#2{(#1)/#2}
\def\afrac#1#2{#1/(#2)}

\newcommand{\eqref}[1]{(\ref{#1})}
\newcommand{\sB}{\mathsf{B}}
\newcommand{\blog}{\mathop{\overline{\log}}}
\newcommand{\KL}{\mathsf{KL}}
\newcommand{\etc}{\textsc{etc}}
\newcommand{\wh}{\widehat}
\newtheorem{lemma}{Lemma}
\newproclaim{definition}{Definition}
\newtheorem{proposition}{Proposition}
\newproclaim{condition}{Condition}
\newtheorem{theorem}{Theorem}

\newproclaim{remark}{Remark}
\newproclaim{exm}{Example}[section]
\newcommand{\1}{\mathbh{1}}
\newcommand{\R}{\mathbb{R}}
\newcommand{\p}{\mathbb{P}}
\newcommand{\E}{\mathbb{E}}
\newcommand{\N}{\mathbb{N}}
\newcommand{\bi}{^{(i)}}
\newcommand{\cI}{\mathcal{I}}
\newcommand{\cJ}{\mathcal{J}}
\newcommand{\cN}{\mathcal{N}}
\newcommand{\cO}{\mathcal{O}}
\newcommand{\cS}{\mathcal{S}}
\newcommand{\cT}{\mathcal{T}}
\newcommand{\cV}{\mathcal{V}}
\newcommand{\eps}{\varepsilon}
\newcommand{\DS}{\displaystyle}
\newlength{\minipagewidth}
\newcommand{\bookbox}[1]{
\par\medskip\noindent
\framebox[\textwidth]{
\begin{minipage}{\minipagewidth}
{#1}
\end{minipage} } \par\medskip}

\newproclaim{ex}{Example}[section]
\newproclaim{exo}{Exercise}[section]
\newtheorem{prop}{Proposition}[section]

\newtheorem{lem}{Lemma}[section]
\newproclaim{rem}{Remark}[section]
\newproclaim{rems}{Remarks}[section]
\newcolumntype{d}{D{.}{.}{6}}
\makeatother

\begin{document}
\begin{frontmatter}

\title{Batched bandit problems}
\runtitle{Batched bandits}

\begin{aug}
% Corresponding author: Philippe Rigollet - rigollet@math.mit.edu% Updated by VTEXPTS2LaTeX.exe, 20.10.2015 08:37
%Updated by VTEXPTS2LaTeX.exe, 19.10.2015 14:31
\author[A]{\fnms{Vianney}~\snm{Perchet}\thanksref{m1,T1}\ead[label=e1]{vianney.perchet@normalesup.org}},
\author[B]{\fnms{Philippe}~\snm{Rigollet}\corref{}\thanksref{m2,T2}\ead[label=e2]{rigollet@math.mit.edu}},
\author[C]{\fnms{Sylvain}~\snm{Chassang}\thanksref{m3,T3}\ead[label=e3]{chassang@princeton.edu}}
\and
\author[D]{\fnms{Erik}~\snm{Snowberg}\thanksref{m4,T3}\ead[label=e4]{snowberg@caltech.edu}}
\runauthor{Perchet, Rigollet, Chassang and Snowberg}
\affiliation{Universit\'e Paris Diderot and INRIA\thanksmark{m1}, Massachusetts
Institute of Technology\thanksmark{m2}, Princeton University\thanksmark{m3} and California Institute
of Technology and NBER\thanksmark{m4}}
%\dedicated{}
\address[A]{V. Perchet\\
LPMA, UMR 7599\\
Universit\'e Paris Diderot\\
8, Place FM/13\\
75013, Paris\\
France\\
\printead{e1}}
\address[B]{P. Rigollet\\
Department of Mathematics and IDSS \\
Massachusetts Institute of Technology\hspace*{27pt}\\
77 Massachusetts Avenue\\
Cambridge, Massachusetts 02139-4307\\
USA\\
\printead{e2}}
\address[C]{S. Chassang\\
Department of Economics\\
Princeton University\\
Bendheim Hall 316\\
Princeton, New Jersey 08544-1021\\
USA\\
\printead{e3}}
\address[D]{E. Snowberg\\
Division of the Humanities and Social Sciences\\
California Institute of Technology\\
MC 228-77\\
Pasadena, California 91125\\
USA\\
\printead{e4}}
\end{aug}
\thankstext{T1}{Supported by ANR Grant ANR-13-JS01-0004.}
\thankstext{T2}{Supported by NSF Grants DMS-13-17308 and CAREER-DMS-10-53987.}
\thankstext{T3}{Supported by NSF Grant SES-1156154.}

% HISTORY:
%
\received{\smonth{5} \syear{2015}}% Updated by VTEXPTS2LaTeX.exe,
%19.10.2015 14:31
%
\revised{\smonth{8} \syear{2015}}% Updated by VTEXPTS2LaTeX.exe,
%19.10.2015 14:31

% ABSTRACT
\begin{abstract}
Motivated by practical applications, chiefly clinical trials, we
study the regret achievable for stochastic bandits under the constraint
that the employed policy must split trials into a small number of
batches. We propose a simple policy, and show that a very small number
of batches gives close to minimax optimal regret bounds. As a
byproduct, we derive optimal policies with low switching cost for
stochastic bandits.
\end{abstract}

% KEYWORDS
\end{frontmatter}

%s1 #&#
\section{Introduction}
\label{SECintro}

All clinical trials are run in \textit{batches}: groups of patients are
treated simultaneously, with the data from each batch influencing the
design of the next. This structure arises as it is impractical to
measure outcomes (rewards) for each patient before deciding what to do
next. Despite the fact that this system is codified into law for drug
approval, it has received scant attention from statisticians. What can
be achieved with a small number of batches? How big should these
batches be? How should results in one batch affect the structure of the next?

We address these questions using the multi-armed bandit framework. This
encapsulates an ``exploration vs. exploitation'' dilemma fundamental to
ethical clinical research \cite{Tho33,Rob52}. In the basic problem,
there are two populations of patients (or \textit{arms}), corresponding to
different treatments. At each point in time $t=1, \ldots, T$, a
decision maker chooses to sample one, and receives a random reward
dictated by the efficacy of the treatment. The objective is to devise a
series of choices---a~policy---maximizing the expected cumulative
reward over $T$ rounds. There is thus a clear tradeoff between
discovering which treatment is the most effective---or \textit{exploration}---and administering the best treatment to as many patients
as possible---or \textit{exploitation}.

The importance of batching extends beyond clinical trials. In recent
years, the bandit framework has been used to study problems in
economics, finance, chemical engineering, scheduling, marketing and,
more recently, internet advertising. This last application has been the
driving force behind a recent surge of interest in many variations of
bandit problems over the past decade. Yet, even in internet
advertising, technical constraints often force data to be considered in
batches; although the size of these batches is usually based on
technical convenience rather than on statistical reasoning. Discovering
the optimal structure, size and number of batches has applications in
marketing \cite{BerMer07,SchBraFad13} and simulations \cite{ChiGan09}.

In clinical trials, batches may be formal---the different phases
required for approval of a new drug by the US~Food and Drug
Administration---or informal---with a pilot, a full trial, and then
diffusion to the full population that may benefit. In an informal
setup, the second step may be skipped if the pilot is successful
enough. In this three-stage approach, the first, and usually second,
phases focus on exploration, while the third focuses on exploitation.
This is in stark contrast to the basic bandit problem described above,
which effectively consists of $T$ batches, each containing a single patient.

We describe a policy that performs well with a small fixed number of
batches. A~fixed number of batches reflects clinical practice, but
presents mathematical challenges. Nonetheless, we identify batch sizes
that lead to a minimax regret bounds as low as the best non-batched
algorithms. We further show that these batch sizes perform well
empirically. Together, these features suggest that near-optimal
policies could be implemented with only small changes to current
clinical practice.

%s2 #&#
\section{Description of the problem}
\label{SECdesc}

%s2.1 #&#
\subsection{Notation}
For any positive integer $n$, define $[n]=\{1,\ldots,n\}$, and for any
$n_1 < n_2$, $[n_1:n_2]=\{n_1, \ldots, n_2\}$ and $(n_1:n_2]=\{n_1+1,
\ldots, n_2\}$. For any positive number $x$, let $\lfloor x \rfloor$
denote the largest integer $n$ such that $n \le x$ and $\lfloor x
\rfloor_2$ denotes the largest \emph{even} integer $m$ such that $m
\le x$. Additionally,\vspace*{1pt} for any real numbers $a$ and $b$, $a\wedge b=\min
(a,b)$ and $a\vee b=\max(a,b)$. Further, define $\blog(x)=1\vee(\log
x)$. $\1(\cdot)$ denotes the indicator function.

If $\cI$, $\cJ$ are closed intervals of $\R$, then $\cI\prec\cJ$
if $x<y$ for all $x \in\cI, y \in\cJ$.

Finally, for two sequences $(u_T)_T, (v_T)_T$, we write $u_T=\cO(v_T)$
or $u_T \lesssim v_T$ if there exists a constant $C>0$ such that
$|u_T|\le C|v_T|$ for any $T$. Moreover, we write $u_T=\Theta(v_T)$ if
$u_T=\cO(v_T)$ and $v_T=\cO(u_T)$.

%s2.2 #&#
\subsection{Framework}
We employ a two-armed bandit framework with horizon $T\ge2$. Central
ideas and intuitions are well captured by this concise framework.
Extensions to $K$-armed bandit problems are mostly technical (see, for
instance, \cite{PerRig13}).

At each\vspace*{1pt} time $t \in[T]$, the decision maker chooses an arm $i \in\{1,
2\}$ and observes a reward that comes from a sequence of i.i.d.~draws
$Y\bi_1, Y\bi_2, \ldots$ from some unknown distribution $\nu\bi$
with expected value $\mu\bi$. We assume that the distributions $\nu
\bi$ are standardized sub-Gaussian, that is, $\int e^{\lambda(x-\mu
\bi)}\nu_i(dx)\le e^{\lambda^2/2}$ for all $\lambda\in\R$. Note
that these include Gaussian distributions with variance at most $1$,
and distributions supported on an interval of length at most $2$.
Rescaling extends the framework to other variance parameters $\sigma^2$.

For any integer $M\in[2:T]$, let $\cT=\{t_1, \ldots, t_{M}\}$ be an
ordered sequence, or \emph{grid}, of integers such that $1<t_1 <\cdots
<t_M=T$. It defines a partition $\cS=\{S_1, \ldots, S_M\}$ of $[T]$
where $S_1=[1:t_1]$ and $S_k=(t_{k-1}:t_{k}]$ for $k \in[2:M]$. The
set $S_k$ is called \emph{$k$th batch}. An \emph{$M$-batch policy} is
a couple $(\cT, \pi)$ where $\cT=\{t_1, \ldots, t_M\}$ is a grid and
$\pi=\{\pi_t , t=1, \ldots, T\}$ is a sequence of random variables
$\pi_t \in\{1, 2\}$, indicating which arm to pull at each time $t=1,
\ldots, T$, which depend only on observations from batches strictly
prior to the current one. Formally, for each $t\in[T]$, let $J(t) \in
[M]$ be the index of the \emph{current batch} $S_{J(t)}$. Then, for $t
\in S_{J(t)}$, $\pi_t$ can only depend on observations $\{Y_s^{(\pi
_s)} :  s \in S_1\cup\cdots\cup S_{J(t)-1}\}=\{Y_s^{(\pi_s)} :
s \le t_{J(t)-1}\}$.

Denote\vspace*{1pt} by $\star\in\{1, 2\}$ the optimal arm defined by $\mu^{(\star
)} = \max_{i \in\{1, 2\}} \mu\bi$, by $\dagger\in\{1, 2\}$ the
suboptimal arm, and by $\Delta:= \mu^{(\star)} - \mu^{(\dagger
)}>0$ the gap between the optimal expected reward and the suboptimal
expected reward.

The performance of a policy $\pi$ is measured by its (cumulative)
\emph{regret} at time $T$
\[
R_T=R_T(\pi)= T \mu^{(\star)} - \sum
_{t=1}^T \E \mu^{({\pi
_t})}.
\]
Denoting by $T_i(t) = \sum_{s=1}^{t} \1(\pi_s = i) , i \in\{1, 2\}$
the number of times arm $i$ was pulled before time $t\geq2$, regret
can be rewritten as $R_T = \Delta\E T_\dagger(T)$.

%s2.3 #&#
\subsection{Previous results}
Bandit problems are well understood in the case where $M=T$, that is,
when the decision maker can use all available data at each time $t \in
[T]$. Bounds on the cumulative regret $R_T$ for stochastic multi-armed
bandits come in two flavors: \emph{minimax} or \emph{adaptive}.
Minimax bounds hold uniformly in $\Delta$ over a suitable subset of
the positive real line such as the intervals $(0,1)$ or even $(0,\infty
)$. The first results of this kind are attributed to Vogel \cite
{Vog60,Vog60b}, who proved that $R_T=\Theta(\sqrt{T})$ in the
two-armed case (see also \cite{FabZwe70,Bat81}).

Adaptive\vspace*{1pt} policies exhibit regret bounds that may be much smaller than
the order of $\sqrt{T}$ when $\Delta$ is large. Such bounds were
proved in the seminal paper of Lai and Robbins \cite{LaiRob85} in an
asymptotic framework (see also \cite{CapGarMai13}). While leading to
tight constants, this framework washes out the correct dependency on
$\Delta$ of the logarithmic terms. In fact, recent research \cite
{AueCesFis02,AudBub10,AueOrt10,PerRig13} has revealed that
$R_T=\Theta (\Delta T\wedge\blog(T\Delta^2)/\Delta )$.

Nonetheless, a systematic analysis of the batched case does not exist,
even though {\sc Ucb2} \cite{AueCesFis02} and {\sc improved-Ucb} \cite
{AueOrt10} are implicitly $M$-batch policies with $M =\Theta(\log T)$.
These algorithms achieve optimal adaptive bounds. Thus, employing a
batched policy is only a constraint when the number of batches $M$ is
much smaller than $ \log T$, as is often the case in clinical practice.
Similarly, in the minimax framework, $M$-batch policies, with $M=\Theta
(\log\log T)$, lead to the optimal regret bound (up to logarithmic
terms) of $\cO (\sqrt{T \log\log\log T} )$ \cite
{CesDekSha13,CesGenMan13}. The sub-logarithmic range $M \ll\log T$ is
essential in applications where $M$ is small and constant, like
clinical trials. In particular, we wish to bound the regret for small
values of $M$, such as 2, 3 or 4.

%s2.4 #&#
\subsection{Literature}
This paper connects to two lines of work: batched sequential estimation
\cite{Dan40,CotJohWet07,Ste45,GhuRob54} and multistage clinical
trials. Somerville \cite{Som54} and Maurice \cite{Mau57} studied the
two-batch bandit problem in a minimax framework under a Gaussian
assumption. They\vspace*{1pt} prove that an ``explore-then-commit'' type policy has
regret of order $T^{2/3}$ for any value of the gap $\Delta$; a result
we recover and extend (see Section~\ref{subminimax}).

Colton \cite{Col63,Col65} introduced a Bayesian perspective,
initiating a long line of work (see \cite{HarSto02} for a recent
overview). Most of this work focuses on the case of two-three batches,
with isolated exceptions \cite{Che96,HarSto02}. Typically,\vspace*{1pt} this work
claims the size of the first batch should be of order $\sqrt{T}$,
which agrees with our results, up to a logarithmic term (see
Section~\ref{subgeom}).

Batched procedures have a long history in clinical trials (see, for instance,
\cite{JenTur00} and \cite{BarLaiShi13}). Usually, batches are of the
same size, or of random size, with the latter case providing
robustness. This literature also focuses on inference questions rather
than cumulative regret. A notable\vspace*{1pt} exception provides an ad-hoc
objective to optimize batch size but recovers the suboptimal $\sqrt
{T}$ in the case of two batches \cite{Bar07}.

%s2.5 #&#
\subsection{Outline}
Section~\ref{SECETC} introduces a general class of $M$-batch policies
we call \textsl{explore-then-commit \textup{(}\etc\textup{)} policies}. These policies
are close to clinical practice within batches. The performance of
generic \etc \ policies are detailed in Proposition~\ref
{PRRegretETC}, found in Section~\ref{SecETCReg}. In Section~\ref{SecFuncGrid}, we study several instantiations of this generic policy
and provide regret bounds with explicit, and often drastic, dependency
on the number of batches $M$. Indeed, in Section~\ref{subminimax},\vadjust{\goodbreak}
we describe a policy in which regret decreases doubly exponentially
fast with the number of batches.

Two\vspace*{1pt} of the instantiations provide adaptive and minimax types of bounds,
respectively. Specifically, we describe two $M$-batch policies, $\pi
^1$ and $\pi^2$ that enjoy the following bounds on the regret:
\begin{eqnarray*}
R_T\bigl(\pi^1\bigr)&\lesssim & \biggl(\frac{T}{\log(T)}
\biggr)^{{1}/{M}}\frac{\blog(T\Delta^2)}{\Delta},
\\
R_T\bigl(\pi^2\bigr)&\lesssim & T^{{1}/({2-2^{1-M}})}
\log^{\alpha_M} \bigl(T^{{1}/({2^M-1})} \bigr), \qquad\alpha_M
\in[0, 1/4).
\end{eqnarray*}
Note that the bound for $\pi^1$ corresponds to the optimal adaptive
rate\break $\blog(T\Delta^2)/\Delta$ when $M=\Theta(\log(T/\log(T)))$
and the bound\vspace*{1.5pt} for $\pi^2$ corresponds to the optimal minimax rate
$\sqrt{T}$ when $M=\Theta(\log\log T)$. The latter is entirely
feasible in clinical settings. As a byproduct of our results, we show
that the adaptive optimal bounds can be obtained with a policy that
switches between arms less than $\Theta(\log(T/\log(T)))$ times,
while the optimal minimax bounds only require $\Theta(\log\log T)$
switches. Indeed, \etc \ policies can be adapted to switch at most once
in each batch.

Section~\ref{SECLB} then examines the lower bounds on regret of any
$M$-batch policy, and shows that the policies identified are optimal,
up to logarithmic terms, within the class of $M$-batch policies.
Finally, in Section~\ref{SECsimulations} we compare policies through
simulations using both standard distributions and real data from a
clinical trial, and show that the policies we identify perform well
even with a very small number of batches.

%s3 #&#
\section{Explore-then-commit policies}
\label{SECETC}

In this section, we describe a simple structure that can be used to
build policies: \emph{explore-then-commit} (\etc). This structure
consists of pulling each arm the same number of times in each
non-terminal batch, and checking after each batch whether, according to
some statistical test, one arm dominates the other. If one dominates,
then only that arm is pulled until $T$. If, at the beginning of the
terminal batch, neither arm has been declared dominant, then the policy
commits to the arm with the largest average past reward. This ``go for
broke'' step is dictated by regret minimization: in the last batch
exploration is pointless as the information it produces can never be used.

Any policy built using this principle is completely characterized by
two elements: the testing criterion and the sizes of the batches.

%s3.1 #&#
\subsection{Statistical test}

We begin by describing the statistical test employed before
non-terminal batches. Denote by
\[
\wh{\mu}\bi_{s} = \frac{1}{s} \sum
_{\ell=1}^s Y\bi_{\ell}
\]
the empirical mean after $s\ge1$ pulls of arm $i$. This estimator
allows for the construction of a collection of upper and lower
confidence bounds\vadjust{\goodbreak} for $\mu\bi$ of the form
\[
\wh{\mu}\bi_{s}+ \sB_{s}\bi \quad\mbox{and} \quad\wh{\mu
}\bi_{s}- \sB_{s}\bi,
\]
where $\sB_s\bi=2\sqrt{2\log (T/s )/s}$ (with the
convention that $\sB_0\bi=\infty$). It follows from Lemma~\ref
{LEMconfbd} that for any $\tau\in[T]$,
%
%e1 #&#
\begin{equation}
\label{EQdeltaCI} \p \bigl\{\exists s\le\tau : \mu\bi> \wh{\mu}\bi_{s}+ \sB
_{s}\bi \bigr\} \vee\p \bigl\{\exists s\le\tau : \mu\bi< \wh {\mu}
\bi_{s}-\sB_{s}\bi \bigr\} \le\frac{4 \tau}{T}.
\end{equation}

These bounds enable us to design the following family of tests $\{
\varphi_t\}_{t \in[T]}$ with values in $\{1,2, \bot\}$ where $\bot$
indicates that the test was inconclusive. This test is only implemented
at times $t \in[T]$ at which each arm has been pulled exactly $s=t/2$
times. However, for completeness, we define the test at all times $t$.
For $t \ge1$, define
\[
\varphi_{t}=\cases{ %
i\in\{1,2\}, & \quad$\!\mbox{if }
T_1(t)=T_2(t)=t/2
\mbox{ and } \wh{\mu}
\bi_{t/2}- \sB_{t/2}\bi> \wh{\mu }^{(j)}_{t/2}+
\sB_{t/2}^{(j)}, j\neq i,\!$ \vspace*{3pt}
\cr
\bot, & \quad$\!
\mbox{otherwise}$.}
\]
The errors\vadjust{\goodbreak} of such tests are controlled as follows.
%
%le3.1 #&#
\begin{lemma}\label{LmConsCrit}
Let $\cS\subset[T]$ be a deterministic subset of even times such that
$T_1(t)=T_2(t)=t/2$, for $t \in\cS$. Partition $\cS$ into $\cS
_-\cup\cS_+$, $\cS_- \prec\cS_+$, where
\[
\cS_-= \biggl\{t \in\cS : \Delta< 16\sqrt{\frac{\log
(2T/t)}{t}} \biggr\},\qquad
\cS_+= \biggl\{t \in\cS : \Delta\ge 16\sqrt{\frac{\log(2T/ t)}{t}} \biggr\} .
\]
Let $\bar t$ denote the smallest element of $\cS_+$. Then
\[
\mathrm{(i)} \ \p(\varphi_{\bar t}\neq\star)\le\frac{4 \bar t}{T} \quad
\mathrm{and} \quad \mathrm{(ii)} \ \p(\exists t \in\cS_- : \varphi_t =
\dagger) \le\frac{4 \bar t}{T}.
\]
\end{lemma}
\begin{pf}
Assume without loss of generality that $\star=1$.
\begin{longlist}[(ii)]
\item[(i)]  By definition,
\[
\{\varphi_{\bar t}\neq1\}= \bigl\{\wh{\mu}^{(1)}_{\bar t/2}-
\sB _{\bar t/2}^{(1)}\le\wh{\mu}^{(2)}_{\bar t/2}+
\sB_{\bar
t/2}^{(2)} \bigr\} \subset\bigl\{ E^1_{\bar t}
\cup E^2_{\bar t} \cup E^3_{\bar t}\bigr\},
\]
where $E^1_t =  \{\mu^{(1)} \geq\wh{\mu}^{(1)}_{t/2}+ \sB
_{t/2}^{(1)} \}$, $E^2_t =  \{\mu^{(2)} \leq\wh{\mu
}^{(2)}_{t/2}- \sB_{t/2}^{(2)}  \}$, and $E^3_t=  \{\mu
^{(1)}-\mu^{(2)}< 2\sB_{t/2}^{(1)}+2\sB_{t/2}^{(2)} \}$. It
follows from~\eqref{EQdeltaCI} that with $\tau=\bar t/2$, $\p
(E_{\bar t}^1)\vee\p(E_{\bar t}^2)\le2 \bar t/T$.

Finally, for any $t \in\cS_+$, in particular for $t=\bar t$, we have
\[
E_t^3 \subset \biggl\{\mu^{(1)}-
\mu^{(2)}< 16\sqrt{\frac{\log(2T/
t)}{t}} \biggr\} =\varnothing.
\]

\item[(ii)]  Focus on the case $t \in\cS_-$, where $ \Delta<
16\sqrt{\log(2T/ t)/t}$. Here,
\[
\bigcup_{t \in\cS_-}\{\varphi_t=2\}=\bigcup
_{t \in\cS_-} \bigl\{ \wh{\mu}^{(2)}_{t/2}-
\sB_{t/2}^{(2)}> \wh{\mu}^{(1)}_{t/2}+ \sB
_{t/2}^{(1)} \bigr\} \subset\bigcup
_{t \in\cS_-}\bigl\{E_t^1 \cup
E_t^2 \cup F_t^3\bigr\},\vadjust{\goodbreak}
\]
\end{longlist}
\noindent where, $E_t^1, E_t^2$ are defined above and $F_t^3= \{\mu^{(1)}-\mu
^{(2)}< 0\}=\varnothing$ as $\star=1$. It follows from~\eqref
{EQdeltaCI}, that with $\tau=\bar t$
\[
\p \biggl( \bigcup_{t \in\cS_-}E_t^1
\biggr)\vee\p \biggl( \bigcup_{t
\in\cS_-}E_t^2
\biggr)\le\frac{2\bar t}{T} .
\]
\upqed\end{pf}

%s3.2 #&#
\subsection{Go for broke}

In the last batch, the \etc \ structure will ``go for broke'' by
selecting the arm $i$ with the largest average. Formally, at time $t$,
let $\psi_t=i$ iff $\wh{\mu}\bi_{T_i(t)}\ge\wh{\mu
}^{(j)}_{T_j(t)}$, with ties broken arbitrarily. While this criterion
may select the suboptimal arm with higher probability than the
statistical test described in the previous subsection, it also
increases the probability of selecting the correct arm by eliminating
inconclusive results. This statement is formalized in the following
lemma. The proof follows immediately from Lemma~\ref{LEMconfbd}.
%
%le3.2 #&#
\begin{lemma}\label{LmIncautCrit}
Fix an even time $t \in[T]$, and assume that both arms have been
pulled $t/2$ times each (i.e., $T_i(t)=t/2$, for $i=1,2$). Going for
broke leads to a probability of error
\[
\p (\psi_t\neq\star )\le\exp\bigl(-t\Delta^2/16\bigr) .
\]
\end{lemma}

%s3.3 #&#
\subsection{Explore-then-commit policy}\label{SecETCReg}

In a batched process, an extra constraint is that past observations can
only be inspected at a specific set of times $\cT=\{t_1, \ldots,
t_{M-1}\}\subset[T]$, called a \emph{grid}.

\begin{figure}
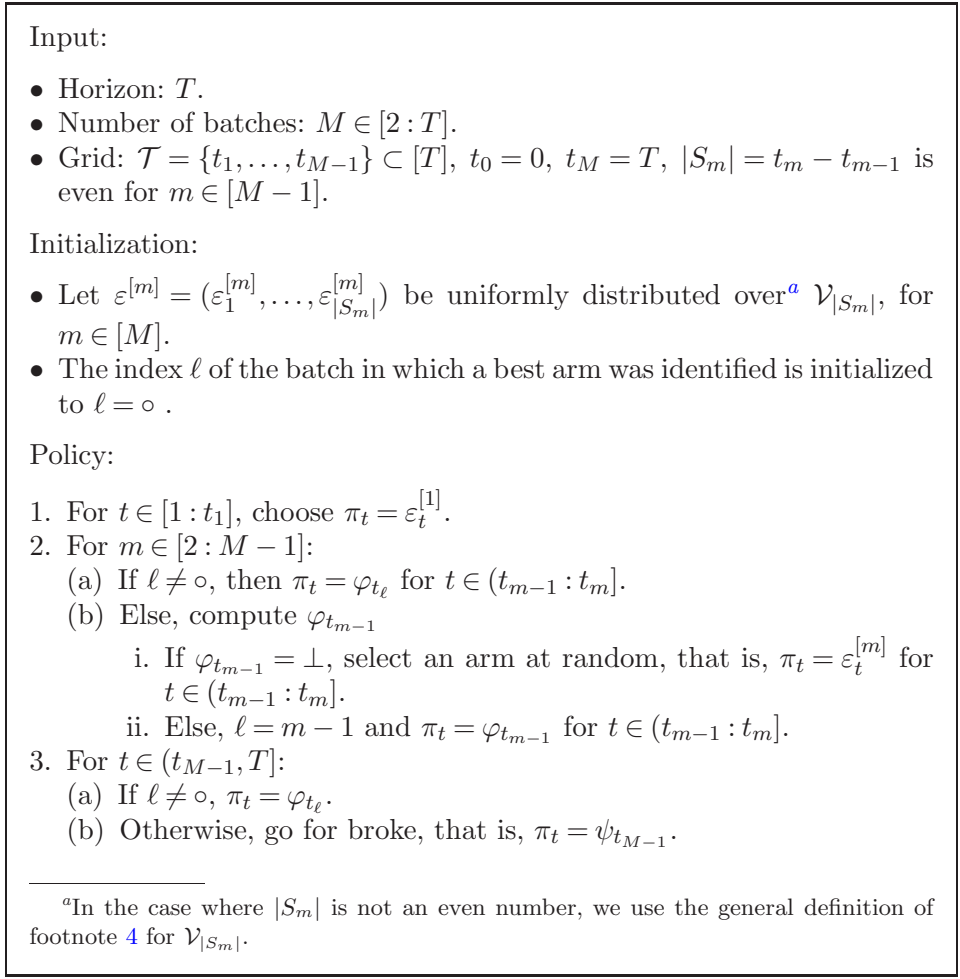

\bookbox{
Input:
\begin{itemize}
\item Horizon: $T$.
\item Number of batches: $M \in[2:T]$.
\item Grid: $\cT=\{t_1, \ldots, t_{M-1}\} \subset[T]$, $t_0=0$,
$t_M=T$, $|S_m|=t_m-t_{m-1}$ is even for $m \in[M-1]$.
\end{itemize}
\noindent Initialization:
\begin{itemize}
\item Let $\eps^{[m]}=(\eps_1^{[m]}, \ldots, \eps_{|S_m|}^{[m]})$
be uniformly distributed over\footnote{In the case where $|S_m|$ is
not an even number, we use the general definition of footnote~\ref
{footunif} for $\cV_{|S_m|}$.} $\cV_{|S_m|}$, for $m \in[M]$.
\item The index $\ell$ of the batch in which a best arm was identified
is initialized to $\ell=\circ$ .
\end{itemize}
Policy:
\begin{enumerate}
\item For $t\in[1:t_1]$, choose $\pi_t=\eps^{[1]}_t$.
\item For $m\in[2: M-1]$:
\begin{enumerate}[(a)]
\item[(a)] If $\ell\neq\circ$, then $\pi_t = \varphi_{t_{\ell}}$ for
$t \in(t_{m-1}:t_{m}]$.
\item[(b)] Else, compute $\varphi_{t_{m-1}}$
\begin{enumerate}
\item If $\varphi_{t_{m-1}}=\bot$, select an arm at random, that is,
$\pi_t=\eps_t^{[m]}$ for $t \in(t_{m-1}:t_{m}]$.
\item Else, $\ell=m-1$ and $\pi_t=\varphi_{t_{m-1}}$ for $t \in
(t_{m-1}:t_{m}]$.
\end{enumerate}
\end{enumerate}
\item For $t \in(t_{M-1},T]$:
\begin{enumerate}[(a)]
\item[(a)] If $\ell\neq\circ$, $\pi_t=\varphi_{t_\ell}$.
\item[(b)] Otherwise, go for broke, that is, $\pi_t=\psi_{t_{M-1}}$.
\end{enumerate}
\end{enumerate}}\vspace*{-6pt}
\caption{Generic explore-then-commit policy with grid $\cT$.}\label{FIGgenericETC}
\end{figure}

The generic \etc \ policy uses a deterministic grid $\cT$ that is
fixed beforehand, and is described more formally in Figure~\ref{FIGgenericETC}. Informally, at each decision time $t_1,\ldots,
t_{M-2}$, the policy implements the statistical test. If one arm is
determined to be better than the other, it is pulled until $T$. If no
arm is declared best, then both arms are pulled the same number of
times in the next batch.
%f1
%f1 #&#
\setcounter{footnote}{3}

We denote by $\eps_t \in\{1,2\}$ the arm pulled at time $t \in[T]$,
and employ an external source of randomness to generate the variables
$\eps_t$. With $N$ an even number, let $(\eps_1, \ldots, \eps_N)$
be uniformly distributed over the subset $\cV_N= \{v \in\{1,2\}
^N :  \sum_{i}\1(v_i=1)=N/2 \}$.\footnote{\label{footunif}Odd
numbers for the deadlines $t_i$ could be considered, at the cost of
rounding problems and complexity, by defining $\cV_N= \{v \in\{
1,2\}^N :   |\sum_{i}\1(v_i=1)-\sum_{i}\1(v_i=2) |\le1
\}$.} This randomization has no effect on the policy, and could easily
be replaced by any other mechanism that pulls each arm an equal number
of times. For example, a mechanism that pulls one arm for the first
half of the batch, and the other for the second half, may be used if
switching costs are a concern.

In the terminal batch $S_M$, if no arm was determined to be optimal in
any prior batch, the \etc \ policy will go for broke by selecting the
arm $i$ such that $\wh{\mu}\bi_{T_i(t_{M-1})}\ge\wh{\mu
}^{(j)}_{T_j(t_{M-1})}$, with\vspace*{1.5pt} ties broken arbitrarily.

To describe the regret incurred by a generic \etc \ policy, we
introduce extra notation. For any $\Delta\in(0,1)$, let $\tau(\Delta
)=T\wedge\vartheta(\Delta)$ where $\vartheta(\Delta)$ is the
smallest integer such that
\[
\Delta\ge16\sqrt{\frac{\log[2T/\vartheta(\Delta)]}{\vartheta
(\Delta)}} .\vadjust{\goodbreak}
\]
Notice that the above definition implies that $\tau(\Delta) \ge2$ and
%
%e2 #&#
\begin{equation}
\label{EQboundNDelta} \tau(\Delta) \le\frac{256}{\Delta^2}\blog \biggl(\frac{T\Delta
^2}{128}
\biggr) .
\end{equation}
The time $\tau(\Delta)$ is, up to a multiplicative constant, the
theoretical time at which the optimal arm will be declared better by
the statistical test with large enough probability. As $\Delta$ is
unknown, the grid will not usually contain this value. Thus, the
relevant time is the first posterior to $\tau(\Delta)$ in a grid:
%
%e3 #&#
\begin{equation}
\label{EQNplus}\quad m(\Delta,\cT) = \cases{ \min\bigl\{m \in\{1,\ldots,M-1\} :
t_{m} \ge \tau(\Delta)\bigr\}, & \quad $\mbox{if } \tau(\Delta) \leq
t_{M-1}$,\vspace*{3pt}
\cr
M-1, & \quad $\mbox{otherwise}$.}\hspace*{-12pt}
\end{equation}

The first proposition gives an upper bound for the regret incurred by a
generic \etc \ policy run with a given set of times $\cT=\{t_1,\ldots
,t_{M-1}\}$.

%pr1 #&#
\begin{proposition}\label{PRRegretETC}
Given the time horizon $T \in\N$, the number of batches $M \in
[2,T]$, and the grid $\cT=\{t_1,\ldots,t_{M-1}\} \subset[T]$ with
$t_0=0$. For any $\Delta\in[0,1]$, the generic \etc \ policy
described in Figure~\ref{FIGgenericETC} incurs regret bounded
%
%e4 #&#
\begin{equation}
\label{EQRegretETC} {R}_T(\Delta, \cT) \leq9\Delta t_{m(\Delta,\cT)} + T
\Delta e^{-({t_{M-1}\Delta^2})/{16}} \1\bigl(m(\Delta,\cT)=M-1\bigr).
\end{equation}
\end{proposition}

\begin{pf}
Denote $\bar m=m(\Delta,\cT)$. Note that $t_{\bar m}$ denotes the
theoretical time on the grid at which the statistical test will declare
$\star$ to be (with high probability) the better arm.

We first examine the case where $t_{\bar m} < M-1$. Define the
following events:
\[
A_m=\bigcap_{n=1}^{m}\{
\varphi_{t_n}=\bot\} , \qquad B_m=\{\varphi
_{t_m}= \dagger\} \quad \mbox{and} \quad C_m=\{
\varphi_{t_{m}}\neq \star\}.
\]
Regret can be incurred in one of the following three manners:
\begin{longlist}[(iii)]
\item[(i)] by exploring before time $t_{\bar m}$,
\item[(ii)] by choosing arm $\dagger$ before time $t_{\bar m}$: this
happens on event $B_m$,
\item[(iii)] by not committing to the optimal arm $\star$ at the
optimal time $t_{\bar m}$: this happens on event $C_{\bar m}$.
\end{longlist}
Error (i) is unavoidable and may occur with probability close to one.
It corresponds to the exploration part of the policy and leads to an
additional term $t_{\bar m}\Delta/2$ in the regret. An error of the
type (ii) or (iii) can lead to a regret of at most $T\Delta$, so we
need to ensure that they occur with low probability. Therefore, the
regret incurred by the policy is bounded as
%
%e5 #&#
\begin{equation}
\label{EQprpropETC1} {R}_T(\Delta, \cT) \le\frac{t_{\bar m}\Delta}{2} + T\Delta\E
\Biggl[\1\Biggl(\bigcup_{m=1}^{\bar m-1}A_{m-1}
\cap B_m\Biggr)+ \1(B_{\bar m -1} \cap C_{\bar m}) \Biggr],
\end{equation}
with the convention that $A_0$ is the whole probability space.

Next, observe that $\bar m$ is chosen such that
\[
16\sqrt{\frac{\log(2T/t_{\bar m})}{t_{\bar m}}}\le\Delta< 16\sqrt {\frac{\log(2T/t_{\bar m-1})}{t_{\bar m-1}}}.
\]
In particular, $t_{\bar m}$ plays the role of $\bar t$ in Lemma~\ref
{LmConsCrit}. Thus, using part (i) of Lemma~\ref{LmConsCrit},
\[
\p(B_{\bar m -1} \cap C_{\bar m})\le\frac{4 t_{\bar m}}{T} .
\]
Moreover, using part (ii) of the same lemma,
\[
\p \Biggl(\bigcup_{m=1}^{\bar m-1}A_{m-1}
\cap B_m \Biggr) \le\frac{4
t_{\bar m}}{T}.
\]
Together with~\eqref{EQprpropETC1} this implies regret is bounded
by ${R}_T(\Delta, \cT) \le9 \Delta t_{\bar m}$.

In the case where $t_{m(\Delta,\cT)} = M-1$, Lemma~\ref
{LmIncautCrit} shows that the go for broke test errs with probability
at most $\exp(-t_{M-1}\Delta^2/16)$, which gives that
\[
{R}_T(\Delta, \cT) \leq9\Delta t_{m(\Delta,\cT)}+T\Delta
e^{-({t_{M-1}\Delta^2})/{16}},
\]
using the same arguments as before.
\end{pf}
Proposition~\ref{PRRegretETC} helps choose a grid by showing how that
choice reduces to an optimal discretization problem.

%s4 #&#
\section{Functionals, grids and bounds}\label{SecFuncGrid}
The regret bound of Proposition~\ref{PRRegretETC} critically depends
on the choice of the grid $\cT=\{t_1,\ldots,t_{M-1}\} \subset[T]$.
Ideally, we would like to optimize the right-hand side of~\eqref
{EQRegretETC} with respect to the $t_m$s. For a fixed $\Delta$, this
problem is easy, and it is enough to choose $M=2$, $t_1 \simeq\tau
(\Delta)$ to obtain optimal regret bounds of the order $R^*(\Delta
)=\log(T\Delta^2)/\Delta$. For unknown $\Delta$, the problem is not
well defined: as observed by \cite{Col63,Col65}, it consists in
optimizing a function ${R}(\Delta, \cT)$ for all $\Delta$, and there
is no choice that is uniformly better than others. To overcome this
limitation, we minimize pre-specified real-valued functionals of
$R(\cdot, \cT)$. The functionals we focus on are:
\begin{eqnarray*}
\DS F_{\mathsf{xs}} \bigl[R_T(\cdot,\cT) \bigr] &=& \sup
_{\Delta\in
[0,1]}\bigl\{R_T(\Delta,\cT)-CR^*(\Delta)\bigr\},
\qquad C>0 \qquad \mbox{Excess regret},
\\
\DS F_{\mathsf{cr}} \bigl[R_T(\cdot,\cT) \bigr]&=&\sup
_{\Delta\in
[0,1]}\frac{R_T(\Delta,\cT)}{R^*(\Delta)}\qquad \mbox{Competitive ratio},
\\
\DS F_{\mathsf{mx}} \bigl[R_T(\cdot,\cT) \bigr]&=& \sup
_{\Delta\in
[0,1]}R_T(\Delta,\cT)\qquad\mbox{Maximum}.
\end{eqnarray*}
Optimizing different functionals leads to different optimal grids. We
investigate the properties of these functionals and grids in the
rest of this section.\footnote{One could also consider the Bayesian
criterion $F_{\mathsf{by}} [R_T(\cdot,\cT) ]=\int R_T(\Delta
,\cT) \,d \pi(\Delta)$ where $\pi$ is a given prior distribution on
$\Delta$, rather than on the expected rewards as in the traditional
Bayesian bandit literature \cite{BerFri85}. %One can also consider
%combination of the Bayesian criterion with other criteria. For example:
%$$
%\int\frac{R_{T}(\Delta,\cT)}{R^*(\Delta)} d \pi(\Delta) .
%$$
}

%s4.1 #&#
\subsection{Excess regret and the arithmetic grid}

We begin with the simple grid consisting in a uniform discretization of
$[T]$. This is particularly prominent in the group sequential testing
literature \cite{JenTur00}. As we will see, even in a favorable setup,
it yields poor regret bounds.

Assume, for simplicity, that $T=2KM$ for some positive integer $K$, so
that the grid is defined by $t_m = mT/M$. In this case, the right-hand
side of~\eqref{EQRegretETC} is bounded \emph{below} by $\Delta
t_1=\Delta T/M$. For small $M$, this lower bound is linear in $T\Delta$,
which is a trivial bound on regret. To obtain a valid upper bound,
note that
\[
t_{m(\Delta, \cT)}\le\tau(\Delta) + \frac{T}{M} \le\frac
{256}{\Delta^2}\blog
\biggl(\frac{T\Delta^2}{128} \biggr) + \frac
{T}{M}.
\]
Moreover, if $m(\Delta,\cT)={M-1}$ then $\Delta$ is of the order of
$\sqrt{1/T}$, thus, $T\Delta\lesssim1/\Delta$. Together with \eqref
{EQRegretETC}, this yields the following theorem.
%
%th1 #&#
\begin{theorem}
\label{THarithmetic}
The \etc \ policy implemented with the arithmetic grid defined above
ensures that, for any $\Delta\in[0,1]$,
\[
{R}_T(\Delta,\cT) \lesssim \biggl(\frac{1}{\Delta}\blog\bigl(T
\Delta^2\bigr) + \frac{T\Delta}{M} \biggr)\wedge T\Delta.
\]
\end{theorem}
The optimal rate is recovered if $M=T$. However, the arithmetic grid
leads to a bound on the excess regret of the order of $\Delta T$ when
$T$ is large and $M$ constant.

In Section~\ref{SECLB}, the bound of Theorem~\ref{THarithmetic} is
shown to be optimal for excess regret, up to logarithmic factors.
Clearly, this criterion provides little useful guidance on how to
attack the batched bandit problem when $M$ is small.

%s4.2 #&#
\subsection{Competitive ratio and the geometric grid}
\label{subgeom}

The geometric grid is defined as $\cT=\{t_1, \ldots, t_{M-1}\}$,
where $t_m=\lfloor a^m \rfloor_2$, and $a\ge2$ is a parameter to be
chosen later. To bound regret using~\eqref{EQRegretETC}, note that if
$m(\Delta,\cT) \leq M-2$, then
\[
{R}_T(\Delta,\cT) \leq9\Delta a^{m(\Delta,\cT)} \le9a\Delta\tau (
\Delta) \le \frac{2304a}{\Delta}\blog \biggl(\frac{T\Delta
^2}{128} \biggr),
\]
and if $m(\Delta,\cT)=M-1$, then $\tau(\Delta) > t_{M-2}$.
Then,~\eqref{EQRegretETC}, together with Lem\-ma~\ref{LEMgeom} yields
\[
{R}_T(\Delta,\cT) \leq9 \Delta a^{M-1}+T\Delta
e^{-\vfrac
{a^{M-1}\Delta^2}{32}} \leq \frac{2336a}{\Delta}\blog \biggl(\frac
{T\Delta^2}{32} \biggr)
\]
for $a \geq2 (\frac{T}{\log T} )^{1/M}\ge2$. We have
proved the following theorem.
%
%th2 #&#
\begin{theorem} \label{thmgeometric}
The \etc \ policy\vspace*{1pt} implemented with the geometric grid defined above for
the value $a:=2 (\frac{T}{\log T} )^{1/M}$, when $M \leq
\log(T/(\log T))$ ensures that, for any $\Delta\in[0,1]$,
\[
{R}_T(\Delta,\cT) \lesssim \biggl(\frac{T}{\log T}
\biggr)^{{1}/{M}}\frac{\blog (T\Delta^2 )}{\Delta} \wedge T\Delta.
\]
\end{theorem}
For a logarithmic number of batches, $M =\Theta(\log T)$, the
geometric grid leads to the optimal regret bound
\[
{R}_T(\Delta,\cT) \lesssim\frac{\blog (T\Delta^2
)}{\Delta} \wedge T\Delta.
\]
This bound shows that the geometric grid leads to a deterioration of
the regret bound by a factor $(T/\log(T))^{1/M}$, which can be
interpreted as a uniform bound on the competitive ratio. For example,
for $M=2$ and $\Delta=1$, this leads to the $\sqrt{T}$ regret bound
observed in the Bayesian literature, which is also optimal in the
minimax sense. However, this minimax optimal bound is not valid for all
values of~$\Delta$. Indeed, maximizing over $\Delta>0$ yields
\[
\sup_\Delta{R}_T(\cT, \Delta)\lesssim
T^{({M+1})/({2M})}\log ^{({M-1})/({2M})}\bigl(\bigl(T/\log(T)\bigr)^{\sfrac{1}{M}}
\bigr),
\]
which yields the minimax rate $\sqrt{T}$ when $M\ge\log(T/\log
(T))$, as expected from prior results. The decay in $M$ can be made
even faster if one focuses on the maximum risk, by employing our
``minimax grid.''

%s4.3 #&#
\subsection{Maximum risk and the minimax grid}
\label{subminimax}

The objective of this grid is to minimize the maximum risk, and to
recover the classical distribution independent minimax bound in $\sqrt
{T}$. The intuition behind this grid comes from Proposition~\ref
{PRRegretETC}, in which $\Delta t_{m(\Delta, \cT)}$ is the most
important term to control. Consider a grid $\cT=\{t_1, \ldots,
t_{M-1}\}$, where the $t_m$'s are defined recursively as
$t_{m+1}=f(t_m)$ so that, by definition, $t_{m(\Delta, \cT)} \le
f(\tau(\Delta)-1)$. As we minimize the maximum risk, $\Delta f(\tau
(\Delta))$ should be the smallest possible term, and constant with
respect to~$\Delta$. This is ensured by choosing $f(\tau(\Delta
)-1)=a/\Delta$ or, equivalently, by choosing $f(x)=a/\tau^{-1}(x+1)$
for a suitable notion of the inverse. This yields \mbox{$\Delta t_{m(\Delta,
\cT)}\le a$}, so that the parameter $a$ is actually a bound on the
regret. This parameter also has to be large enough so that the regret
$T\sup_{\Delta} \Delta e^{-t_{M-1}\Delta^2/8}=2T/\sqrt{et_{M-1}}$
incurred in the go for broke step is also of the order of $a$. The
formal definition below uses not only this delicate recurrence, but
also takes care of rounding problems.

Let $u_1=a$, for some $a>0$ to be chosen later, and $u_{j}=f(u_{j-1})$ where
%
%e6 #&#
\begin{equation}
\label{EQRecMinMax} f(u)=a\sqrt{\frac{u}{\log (({2T})/{u} )}}
\end{equation}
for all $j \in\{2,\ldots, M-1\}$. The \emph{minimax grid} $\cT=\{
t_1, \ldots, t_{M-1}\}$ has points given by $t_m=\lfloor u_m\rfloor
_2, m\in\{1,\ldots, M-1\}$.

If $m(\Delta, \cT)\le M-2$, then it follows from~\eqref
{EQRegretETC} that $R_T(\Delta, \cT)\le9\Delta t_{m(\Delta, \cT
)}$, and as $\tau(\Delta)$ is the smallest integer such that $\Delta
\ge16a/f(\tau(\Delta))$, we have
\[
\Delta t_{m(\Delta, \cT)}\le\Delta f\bigl(\tau(\Delta)-1\bigr)\le16a .
\]
As discussed above, if $a$ is greater than $2\sqrt{2}T/(16\sqrt
{et_{M-1}})$, then the regret is also bounded by $16a$ when $m(\Delta,
\cT)=M-1$. Therefore, in both cases, the regret is bounded by $16a$.
Before finding\vspace*{1pt} an $a$ satisfying the above conditions, note that it
follows from Lemma~\ref{LEMinduction2} that, as long as
$15a^{S_{M-2}} \leq2T$,
\[
t_{M-1}\ge\frac{u_{M-1}}{2}\ge\frac{a^{S_{M-2}}}{30\log^{{S_{M-3}}/{2}} (2T/a^{S_{M-5}} )} ,
\]
with the notation $S_k:=2-2^{-k}$. Therefore, we need to choose $a$
such that
\[
a^{S_{M-1}}\ge\sqrt{\frac{15}{16e}}T \log^{{S_{M-3}}/{4}} \biggl(
\frac{2T}{a^{S_{M-5}}} \biggr)\quad \mbox{and}\quad 15a^{S_{M-2}} \leq2T .
\]
It follows from Lemma~\ref{LEMChoicea} that the choice
\[
a:=(2T)^{{1}/{S_{M-1}}}\log^{{1}/{4}-({3}/{4}){1}/({2^M-1})} \bigl((2T)^{{15}/({2^M-1})} \bigr)
\]
ensures both conditions when $2^M \leq\log(2T)/6$. We emphasize that
\[
\log^{{1}/{4}-({3}/{4}){1}/({2^M-1})} \bigl((2T)^{{15}/({2^M-1})} \bigr) \leq2 \qquad\mbox{ with }
M=\bigl\lfloor\log _2\bigl(\log(2T)/6\bigr)\bigr\rfloor.
\]
As a consequence, in order to get the optimal minimax rate of $\sqrt
{T}$, one only needs $\lfloor\log_2\log(T) \rfloor$ batches. If
more batches are available, then our policy implicitly combines some of
them. We have proved the following theorem.
%
%th3 #&#
\begin{theorem} \label{thmminimax}
The \etc \  policy over the minimax grid with
\[
a= (2T)^{\afrac{1}{2-2^{1-M}}}\log^{\sfrac{1}{4}-(\sfrac{3}{4})\afrac
{1}{2^M-1}} \bigl((2T)^{\afrac{15}{2^M-1}} \bigr)
\]
ensures that, for any $M$ such that $2^M \leq\log(2T)/6$,
\[
\sup_{0\le\Delta\le1}{R}_T(\Delta,\cT) \lesssim
T^{\afrac
{1}{2-2^{1-M}}}\log^{\sfrac{1}{4}-(\sfrac{3}{4})\afrac{1}{2^M-1}} \bigl(T^{\afrac{1}{2^M-1}} \bigr),
\]
which is minimax optimal, that is, $\sup_\Delta{R}_T(\Delta,\cT)
\lesssim\sqrt{T}$, for $M \geq \log_2\log(T)$.
\end{theorem}

Table~\ref{TABminimax} gives the regret bounds (without constant
factors) and the decision times of the \etc \ policy with the minimax
grid for $M=2,3,4,5$.

%t1 #&#
\begin{table}
\caption{Regret and decision times of the \etc \ policy with the
minimax grid for $M=2,3,4,5$. In the table, $l_T=\log(T)$}
\label{TABminimax}
\begin{tabular*}{\textwidth}{@{\extracolsep{\fill}}lcccc@{}}
\hline
$\bolds{M}$ &$\bolds{t_1=\sup_{\Delta} R_T(\Delta,\cT)}$& $\bolds{t_2}$ &
$\bolds{t_3}$& $\bolds{t_4}$\\
\hline
$2$ & $T^{2/3}$& & &\\[3pt]
$3$ & $T^{4/7}l_T^{1/7}$& $ T^{6/7}l_T^{-1/7}$ & & \\[3pt]
$4$ &$T^{8/15}l_T^{1/5}$& $T^{12/15}l_T^{-1/5}$ & $T^{14/15}l_T^{-2/5}$
&\\[3pt]
$5$ &$T^{16/31}l_T^{7/31}$& $T^{24/31}l_T^{-5/31}$ &
$T^{28/31}l_T^{-11/31}$ &$T^{30/31}l_T^{-14/31}$\\
\hline
\end{tabular*}
\end{table}

The \etc \ policy with the minimax grid can easily be adapted to have
only $O(\log\log T)$ switches, and yet still achieve regret of optimal
order $\sqrt{T}$. To do so, in each batch one arm should be pulled for
the first half of the batch, and the other for the second half, leading
to only one switch within the batch, until the policy commits to a
single arm. To ensure that a switch does not occur between batches, the
first arm pulled in a batch should be set to the last arm pulled in the
previous batch, assuming that the policy has not yet committed. This
strategy is relevant in applications such as labor economics and
industrial policy, where switching from an arm to the other may be
expensive \cite{Jun04}. In this context, our policy compares favorably
with the best current policies constrained to have $\log_2\log(T)$
switches, which lead to a regret bound of order $\sqrt{T \log\log
\log T}$ \cite{CesDekSha13}.

%s5 #&#
\section{Lower bounds}
\label{SECLB}
In this section, we address the optimality of the regret bounds derived
above for the specific functionals $F_{\mathsf{xs}}$, $F_{\mathsf{cr}}$ and $F_{\mathsf{mx}}$. The results below do not merely
characterize optimality (up to logarithmic terms) of the chosen grid
within the class of \etc \ policies, but also optimality of the final
policy among the class of \emph{all $M$-batch policies}.

%th4 #&#
\begin{theorem}\label{THlb}
Fix $T\ge2$ and $M \in[2:T]$. Any $M$-batch policy $(\cT, \pi)$,
must satisfy the following lower bounds:
\begin{eqnarray*}
\DS\sup_{\Delta\in(0,1]} \biggl\{R_T(\Delta, \cT) -
\frac
{1}{\Delta_{\mathsf{xs}}} \biggr\}&\gtrsim & \frac{T}{M},
\\
\DS\sup_{\Delta\in(0,1]} \bigl\{\Delta R_T(\Delta, \cT)
\bigr\} &\gtrsim & T^{\sfrac{1}{M}},
\\
\DS\sup_{\Delta\in(0,1]} \bigl\{R_T(\Delta, \cT) \bigr\} &
\gtrsim & T^{\afrac{1}{2-2^{1-M}}}.
\end{eqnarray*}
\end{theorem}
\begin{pf}
Fix\vspace*{1pt} $\Delta_k=\frac{1}{\sqrt{t_k}}, k=1,\ldots, M$. Focusing first
on excess risk, it follows from Proposition~\ref{proplb} that
\begin{eqnarray*}
&&\sup_{\Delta\in(0,1]} \biggl\{R_T(\Delta, \cT)-
\frac{1}{\Delta} \biggr\}\\
&&\qquad\ge  \max_{1\le k \le M}\sum
_{j=1}^M \biggl\{\frac{\Delta_k
t_j}{4}\exp
\bigl(-t_{j-1}\Delta_k^2/2\bigr) -
\frac{1}{\Delta_k} \biggr\}
\\
&&\qquad\ge  \max_{1\le k \le M} \biggl\{\frac{t_{k+1}}{4\sqrt{et_k}}-\sqrt
{t_k} \biggr\}.
\end{eqnarray*}
As $t_{k+1} \ge t_k$, the last quantity above is minimized if all the
terms are of order 1. This yields $t_{k+1}=t_k + a $, for some
positive constant $a$. As $t_M=T$, we get that $t_j \sim jT/M$, and
taking $\Delta=1$ yields
\[
\sup_{\Delta\in(0,1]} \biggl\{R_T(\Delta, \cT)-
\frac{1}{\Delta} \biggr\}\ge\frac{t_1}{4} \gtrsim\frac{T}{M} .
\]
Proposition~\ref{proplb} also yields
\begin{eqnarray*}
 \sup_{\Delta\in(0,1]} \bigl\{\Delta R_T(\Delta, \cT)
\bigr\}&\ge &\max_{k}\sum_{j=1}^M
\biggl\{\frac{\Delta_k^2 t_j}{4}\exp\biggl(-\frac
{t_{j-1}\Delta_k^2}{2}\biggr) \biggr\}\\
&\ge&\max
_{k} \biggl\{\frac
{t_{k+1}}{4\sqrt{e}t_k} \biggr\}.
\end{eqnarray*}
Arguments similar to the ones for the excess regret above, give the
lower bound for the competitive ratio. Finally,
\begin{eqnarray*}
\sup_{\Delta\in(0,1]}R_T(\Delta, \cT)&\ge&\max
_{k}\sum_{j=1}^M
\biggl\{\frac{\Delta_k t_j}{4}\exp\biggl(-\frac{t_{j-1}\Delta_k^2}{2}\biggr)
\biggr\}\\
&\ge&\max
_{k} \biggl\{\frac{t_{k+1}}{4\sqrt{et_k}} \biggr\}
\end{eqnarray*}
gives the lower bound for maximum risk.
\end{pf}

%s6 #&#
\section{Simulations}
\label{SECsimulations}
In this final section, we briefly compare, in simulations, the various
policies (grids) introduced above. These are also compared with {\sc
Ucb2} \cite{AueCesFis02}, which, as noted above, can be seen as an
$M=O(\log T)$ batch trial. A more complete exploration can be found in
\cite{PerRigCha15supp}.

The minimax and geometric grids perform well using an order of
magnitude fewer batches than {\sc Ucb2}. The number of batches required
for {\sc Ucb2} make its use for medical trials functionally impossible.
For example, a study that examined STI status six months after an
intervention in~\cite{MetFeaGoo13} would require 1.5 years to run
using minimax batch sizes, but {\sc Ucb2} would use as many as 56
batches, meaning the study would take 28 years.

%f2 #&#
\begin{figure}

\includegraphics{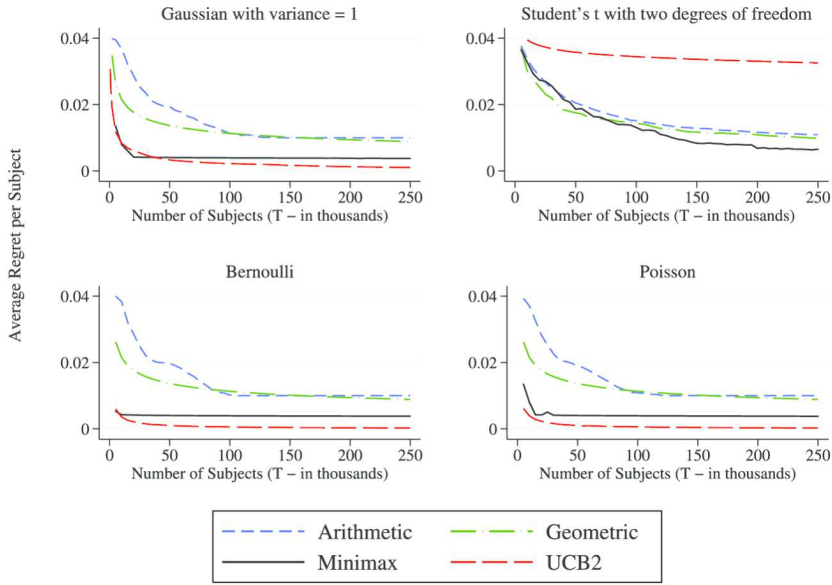}

\caption{Performance of policies with different distributions and
$M=5$. (For all distributions $\mu^{(\dagger)}=0.5$, and $\mu
^{(\star)} = 0.5+\Delta= 0.6$.)}\label{figdistributions}
\end{figure}

Specific examples of performance can be found in Figure~\ref{figdistributions}. This figure compares average regret produced by
different policies and many values of the total sample, $T$. For each
value of $T$ in the figure, a sample is drawn, grids are computed based
on $M$ and $T$, the policy is implemented, and average regret is
calculated based on the choices in the policy. This is repeated 100
times for each value of $T$.

The number of batches is set at $M=5$ for all policies except {\sc
Ucb2}. Each panel considers one of four distributions: two
continuous---Gaussian and Student's $t$-distribution---and two
discrete---Bernoulli and Poisson. In all cases, we set the difference
between the arms at $\Delta= 0.1$.

A few patterns are immediately apparent. First, the arithmetic grid
produces relatively constant average regret above a certain number of
participants. The intuition is straightforward: when $T$ is large
enough, the \etc \ policy will tend to commit after the first batch, as
the first evaluation point will be greater than $\tau(\Delta)$. In
the arithmetic grid, the size of this first batch is a constant
proportion of the overall participant pool, so average regret will be
constant when $T$ is large enough.

Second, the minimax grid also produces relatively constant average
regret, although this holds for smaller values of $T$, and produces
lower regret than the geometric or arithmetic case when $M$ is small.
This indicates, using the intuition above, that the minimax grid excels
at choosing the optimal batch size to allow a decision to commit very
close to $\tau(\Delta)$. This advantage over the arithmetic and
geometric grids is clear. The minimax grid can even produce lower
regret than {\sc Ucb2}, using an order of magnitude fewer batches.

Third, and finally, the {\sc Ucb2} algorithm generally produces lower
regret than any of the policies considered in this manuscript for all
distributions except the heavy-tailed Student's $t$-distribution, for
which batched policies perform significantly better. Indeed, the {\sc
Ucb2} is calibrated for sub-Gaussian rewards, as are batched policies.
However, even with heavy-tailed distributions, the central limit
theorem implies that batching a large number of observations returns
averages that are sub-Gaussian; see the supplementary material \cite
{PerRigCha15supp}. Even when {\sc Ucb2} performes better, this increase
in performance comes at a steep practical cost: many more batches. For
example, with draws from a Gaussian distribution, and $T$ between
10{,}000 and 40{,}000, the minimax grid with only $5$ batches performs
better than {\sc Ucb2}. Throughout this range, {\sc Ucb2} uses roughy
50 batches.

It is worth noting that in medical trials, there is nothing special
about waiting six months for data from an intervention. Trials of
cancer drugs often measure variables like the 1- or 3-year survival
rate, or the increase in average survival compared to a baseline that
may be greater than a year. In these cases, the ability to get
relatively low regret with a small number of batches is extremely important.

%sA #&#
\begin{appendix}
%sB #&#
\section{Tools for lower bounds}

Our results hinge on tools for lower bounds, recently adapted to the
bandit setting in \cite{BubPerRig13}. Specifically, we reduce the
problem of deciding which arm to pull to that of hypothesis testing.
Consider the following two candidate setups for the rewards
distributions: $P_1=\cN(\Delta, 1)\otimes\cN(0,1)$ and $P_2=\cN(0,
1)\otimes\cN(\Delta,1)$, that is, under $P_1$ successive pulls of
arm~1 yield $\cN(\Delta,1)$ rewards and successive pulls of arm~2
yield $\cN(0,1)$ rewards. The opposite is true for $P_2$, so arm $i$
is optimal under $P_i$.

At a given time $t \in[T]$, the choice of $\pi_t \in\{1, 2\}$ is a
test between $P_1^t$ and $P_2^t$ where $P_i^t$ denotes the distribution
of observations available at time $t$ under $P_i$. Let $R(t, \pi)$
denote the regret incurred by policy $\pi$ at time $t$. We have $R(t,
\pi)=\Delta\1(\pi_t\neq i)$. Denote by $E_i^t$ the expectation under
$P_i^t$, so that
\begin{eqnarray*}
E_1^t \bigl[R(t, \pi) \bigr]\vee E_2^t
\bigl[R(t, \pi) \bigr]&\ge& \frac
{1}{2} \bigl(E_1^t
\bigl[R(t, \pi) \bigr]+ E_2^t \bigl[R(t, \pi) \bigr]
\bigr)
\\
&=& \frac{\Delta}{2} \bigl(P_1^t(\pi_t=2)+
P_2^t(\pi_t=1) \bigr).
\end{eqnarray*}

Next, we use the following lemma (see \cite{Tsy09}, Chapter~2).
%
%leB.1 #&#
\begin{lem}
\label{lemlbtest}
Let $P_1$ and $P_2$ be two probability distributions such that \mbox{$P_1 \ll
P_2$}. Then for any measurable set $A$,
\[
P_1(A)+P_2\bigl(A^c\bigr)\ge
\tfrac{1}{2}\exp \bigl(-\KL(P_1,P_2) \bigr),
\]
where $\KL(\cdot,\cdot)$ is the Kullback--Leibler divergence defined by
\[
\KL(P_1,P_2)=\int\log \biggl(\frac{d P_1}{d P_2} \biggr)
\,d P_1.
\]
\end{lem}

Here, observations are generated by an $M$-batch policy $\pi$. Recall
that $J(t) \in[M]$ denotes the index of the current batch. As $\pi$
depends on observations $\{Y_s^{(\pi_s)} :  s\in[t_{J(t)-1}]\}$,
$P_i^t$ is a product distribution of at most $t_{J(t)-1}$ marginals. It
is straightforward to show that whatever arms are observed over the
history, $\KL(P_1^{t},P_2^{t})=t_{J(t)-1}\Delta^2/2$. Therefore,
\[
E_1^t \bigl[R(t, \pi) \bigr]\vee E_2^t
\bigl[R(t, \pi) \bigr] \ge\tfrac
{1}{4}\exp \bigl(-t_{J(t)-1}
\Delta^2/2 \bigr).
\]
Summing over $t$ yields the following result.
%
%prB.1 #&#
\begin{prop}
\label{proplb}
Fix $\cT=\{t_1, \ldots, t_M\}$ and let $(\cT, \pi)$ be an $M$-batch
policy. There exist reward distributions with gap $\Delta$, such that
$(\cT, \pi)$ has regret bounded below as, defining $t_0:=0$,
\[
R_T(\Delta, \cT)\ge\Delta\sum_{j=1}^M
\frac{t_j}{4}\exp \bigl(-t_{j-1}\Delta^2/2 \bigr).
\]
\end{prop}
\noindent A variety of lower bounds in Section~\ref{SECLB} are shown
using this proposition.

%sC #&#
\section{Technical lemmas}

A process $\{Z_t\}_{t\ge0}$ is a sub-Gaussian martingale difference
sequence if $\E [Z_{t+1} | Z_1,\dots,Z_t ]=0$ and $\E
[e^{\lambda Z_{t+1}} ]\leq e^{\lambda^2/2}$ for every $\lambda>0,
t\ge0$.

%leC.1 #&#
\begin{lem}\label{LEMconfbd}
Let $Z_t$ be a sub-Gaussian martingale difference sequence. Then, for
every $\delta>0$ and every integer $t \ge1$,
\[
\p \biggl\{\bar{Z}_t \geq\sqrt{\frac{2}{t} \log \biggl(
\frac
{1}{\delta} \biggr)} \biggr\} \leq\delta.
\]
Moreover, for every integer $\tau\ge1$,
\[
\p \biggl\{\exists t \leq\tau, \bar{Z}_t \geq2\sqrt{
\frac
{2}{t} \log \biggl(\frac{4}{\delta}\frac{\tau}{t} \biggr)}
\biggr\} \leq\delta.
\]
\end{lem}
\begin{pf}
The first inequality follows from a classical Chernoff
bound. To prove the maximal inequality, define $\varepsilon_t=2\sqrt
{\frac{2}{t} \log (\frac{4}{\delta}\frac{\tau}{t} )}$.
Note that, by Jensen's inequality, for any $\alpha>0$, the process $\{
\exp(\alpha s\bar Z_s)\}_s$ is a sub-\break martingale. Therefore, it follows
from Doob's maximal inequality \cite{Doo90}, Theorem~3.2, page~314,  that
for every $\eta>0$ and every integer $t \ge1$,
\begin{eqnarray*}
\p \{\exists s \leq t, s\bar{Z}_s \geq\eta \}&=&\p \bigl\{\exists s
\leq t, e^{\alpha s\bar{Z}_s} \geq e^{\alpha\eta
} \bigr\}\\
& \leq&\E
\bigl[e^{\alpha t\bar Z_t} \bigr]e^{-\alpha\eta
}.
\end{eqnarray*}

Next, as $Z_t$ is sub-Gaussian, we have $\E [\exp(\alpha t\bar
Z_t) ] \le\exp(\alpha^2t/2)$. The above, and optimizing with
respect to $\alpha>0$ yields
\[
\p \{\exists s \leq t, s\bar{Z}_s \geq\eta \} \leq \exp \biggl(-
\frac{\eta^2}{2t} \biggr).
\]
Next, using a peeling argument, one obtains
\begin{eqnarray*}
\p \{\exists t \leq\tau, \bar{Z}_t \geq\varepsilon _t \}
&\leq & \sum_{m=0}^{\lfloor\log_2(\tau)\rfloor} \p \Biggl\{ \bigcup
_{ t =2^m}^{ 2^{m+1}-1}\{ \bar{Z}_t \geq
\varepsilon_t\} \Biggr\}
\\
&\leq& \sum_{m=0}^{\lfloor\log_2(\tau)\rfloor} \p \Biggl\{\bigcup
_{
t =2^m}^{ 2^{m+1}}\{ \bar{Z}_t
\geq
\varepsilon_{2^{m+1}}\} \Biggr\} \\
&\leq & \sum_{m=0}^{\lfloor\log_2(\tau)\rfloor}
\p \Biggl\{\bigcup_{
t =2^m}^{ 2^{m+1}}\bigl\{ t
\bar{Z}_t \geq2^m\varepsilon_{2^{m+1}}\bigr\} \Biggr
\}
\\
& \leq & \sum_{m=0}^{\lfloor\log_2(\tau)\rfloor} \exp \biggl(-
\frac
{ (2^m\varepsilon_{2^{m+1}} )^2}{2^{m+2}} \biggr)\\
&=&\sum_{m=0}^{\lfloor\log_2(\tau)\rfloor}
\frac{2^{m+1}}{\tau}\frac
{\delta}{4}\\
&\le&\frac{2^{\log_2(\tau)+2}}{\tau}\frac{\delta
}{4}\le
\delta.
\end{eqnarray*}
Hence, the result.
\end{pf}

%leC.2 #&#
\begin{lem}
\label{LEMgeom}
Fix two positive integers $T$ and $M \leq\log(T)$. It holds that
\[
T\Delta e^{-\vfrac{a^{M-1}\Delta^2}{32}}\leq32 a \frac{\blog
(\vfrac{T\Delta^2}{32} )}{\Delta} \qquad \mbox{if } a \geq
\biggl(\frac{MT}{\log T} \biggr)^{\sfrac{1}{M}}.
\]
\end{lem}
\begin{pf}
Fix the value of $a$ and observe that $M \leq\log T$
implies that $a \geq e$. Define $x:=T\Delta^2/32 > 0$ and $\theta:=
a^{M-1}/T > 0$. The first inequality is rewritten as
%
%eC.1 #&#
\begin{equation}
\label{EQLemGeom} xe^{- \theta x} \leq a \blog(x).
\end{equation}
We will prove that this inequality is true for all $x>0$, given that
$\theta$ and $a$ satisfy some relation. This, in turn, gives a
condition that depends solely on $a$, ensuring that the statement of
the lemma is true for all $\Delta>0$.

Equation \eqref{EQLemGeom} immediately holds if $x \leq e$ as $a\blog
(x)=a\geq e$. Similarly, $xe^{-\theta x} \leq1/(\theta e)$. Thus
\eqref{EQLemGeom} holds for all $x \geq1/\sqrt{\theta}$ when $a
\geq a^* : =1/(\theta\blog(1/\theta))$. We assume this inequality
holds. Thus, we must show that \eqref{EQLemGeom} holds for $x \in[e,
1/\sqrt{\theta}]$. For $x \leq a$, the derivative of the right-hand
side is $\frac{a}{x} \geq1$, while the derivative of the left-hand
side is smaller than 1. As a consequence, \eqref{EQLemGeom} holds for
every $x \leq a$, in particular for every $x \leq a^*$. To summarize, whenever
\[
a \geq a^*=\frac{T}{a^{M-1}}\frac{1}{\blog(\sfrac{T}{a^{M-1}})},
\]
equation\vspace*{1pt} \eqref{EQLemGeom} holds on $(0,e]$, on $[e,a^*]$ and on $[1/\sqrt
{\theta},+\infty)$, thus on $(0,+\infty)$ as $a^* \geq1/\sqrt
{\theta}$. Next, if $a^M \geq MT/ \log T$, we obtain
\begin{eqnarray*}
\frac{a}{a^*} &=& \frac{a^M}{T} \blog \biggl(\frac{T}{a^{M-1}}
\biggr)\\
&\geq&\frac{M}{\log(T)}\log \biggl( T \biggl(\frac{\log T}{M
T}
\biggr)^{\vfrac{M-1}{M}} \biggr)
\\
&=&\frac{1}{\log(T)} \log \biggl(T \biggl(\frac{\log(T)}{M}
\biggr)^{M-1} \biggr).
\end{eqnarray*}
The result follows from $\log(T)/M \geq1$, hence $a/a^* \geq1$.
\end{pf}

%leC.3 #&#
\begin{lem}\label{LEMinduction2}
Fix $a \ge1, b \ge e$ and let $u_1,u_2, \ldots$ be defined by $u_1=
a$ and $u_{k+1}=a\sqrt{\frac{u_k}{\log(b/u_k)}}$. Define $S_k=0$ for
$k <0$ and
\[
S_k=\sum_{j=0}^k2^{-j}=2-2^{-k}
\qquad\mbox{for }k \ge0.
\]
Then, for any $M$ such that $15a^{S_{M-2}}\le b$, and all $k \in[M-3]$,
\[
u_k \ge \frac{a^{S_{k-1}}}{15\log^{\sfrac{S_{k-2}}{2}}
(b/a^{S_{k-2}} )}.
\]
Moreover, for $k \in[M-2:M]$, we also have
\[
u_k \ge \frac{a^{S_{k-1}}}{15\log^{\sfrac{S_{k-2}}{2}}
(b/a^{S_{M-5}} )}.
\]
\end{lem}
\begin{pf}
Define\vspace*{1pt} $z_k=\log(b/a^{S_{k}})$. It is straightforward to show that
$z_k \le3z_{k+1}$ iff $a^{S_{k+2}} \le b$. In particular,
$a^{S_{M-2}}\le b$ implies that $z_k \le3 z_{k+1}$ for all $k \in
[0: M-4]$. Next, we have
%
%eC.2 #&#
\begin{equation}
\label{EQukplus1} u_{k+1}=a\sqrt{\frac{u_k}{\log(b/u_k)}}\ge a\sqrt{
\frac {a^{S_{k-1}}} {15z_{k-2}^{\sfrac{S_{k-2}}{2}}\log(b/u_k)}}.
\end{equation}
Observe that $b/a^{S_{k-1}}\ge15$, so for all $k \in[0,M-1]$ we have
\begin{eqnarray*}
&&\log(b/u_k)\le\log\bigl(b/a^{S_{k-1}}\bigr) + \log15+
\frac{S_{k-2}}{2} \log z_{k-2} \le5z_{k-1}.
\end{eqnarray*}
This yields
\[
z_{k-2}^{\sfrac{S_{k-2}}{2}}\log(b/u_k) \le15z_{k-1}^{\sfrac
{S_{k-2}}{2}}z_{k-1}=15z_{k-1}^{S_{k-1}}.
\]
Plugging this bound into~\eqref{EQukplus1} completes the proof for
$k\in[M-3]$.

Finally, if $k \ge M-2$, we have by induction on $k$ from $M-3$,
\[
u_{k+1}=a\sqrt{\frac{u_k}{\log(b/u_k)}}\ge a\sqrt{
\frac{a^{S_{k-1}}}{15z_{M-5}^{\sfrac{S_{k-2}}{2}}\log(b/u_k)}}.
\]
Moreover, as $b/a^{S_{k-1}} \ge15$, for $k \in[M-3, M-1]$ we have
\[
\log(b/u_k)\le\log\bigl(b/a^{S_{k-1}}\bigr) + \log15+
\frac{S_{k-2}}{2} \log z_{M-5} \le3z_{M-5}.
\]
\upqed\end{pf}
%
%leC.4 #&#
\begin{lem}\label{LEMChoicea}
If $2^{M} \leq\log(4T)/6$, the
following specific choice
\[
a:=(2T)^{\sfrac{1}{S_{M-1}}}\log^{\sfrac{1}{4}-(\sfrac{3}{4})\afrac
{1}{2^M-1}} \bigl((2T)^{\afrac{15}{2^M-1}} \bigr)
\]
ensures that
%
%eC.3 #&#
\begin{equation}
\label{EQa41} a^{S_{M-1}}\ge\sqrt{\frac{15}{16e}}T
\log^{\sfrac{S_{M-3}}{4}} \biggl(\frac{2T}{a^{S_{M-5}}} \biggr)
\end{equation}
and
%
%eC.4 #&#
\begin{equation}
\label{EQa42} 15a^{S_{M-2}} \leq2 T.
\end{equation}
\end{lem}
\begin{pf}
Immediate for $M=2$. For $M >2$, $2^{M} \leq\log(4T)$ implies
\begin{eqnarray*}
a^{S_{M-1}} &=&2T \log^{\sfrac{S_{M-3}}{4}} \bigl((2T)^{\afrac
{15}{2^M-1}} \bigr)
\\
&\geq& 2T \biggl[16\frac{15}{2^M-1} \log(2T) \biggr]^{1/4}\ge2T.
\end{eqnarray*}
Therefore, $a \geq(2T)^{1/S_{M-1}}$, which in turn implies that
\begin{eqnarray*}
a^{S_{M-1}} &=& 2T \log^{\sfrac{S_{M-3}}{4}} \bigl((2T)^{1-\sfrac
{S_{M-5}}{S_{M-1}}} \bigr)\\
&\geq&
\sqrt{\frac{15}{16e}}T \log^{\sfrac
{S_{M-3}}{4}} \biggl(\frac{2T}{a^{S_{M-5}}} \biggr).
\end{eqnarray*}
This completes the proof of~\eqref{EQa41}. Equation \eqref{EQa42} follows if
%
%eC.5 #&#
\begin{equation}
\label{EQa43} 15^{S_{M-1}}(2T)^{S_{M-2}}\log^{\vfrac{S_{M-3}S_{M-2}}{4}}
\bigl((2T)^{\afrac{15}{2^M-1}} \bigr) \leq(2T)^{S_{M-1}}.
\end{equation}
Using that $S_{M-k} \leq2$, we get that the left-hand side of~\eqref
{EQa42} is smaller than
\[
15^2 \log \bigl((2T)^{\afrac{15}{2^M-1}} \bigr)\leq2250\log
\bigl((2T)^{2^{1-M}} \bigr).
\]
The result follows using $2^M \leq\log(2T)/6$, which implies that the
right-hand side in the above inequality is bounded by $(2T)^{2^{1-M}}$.
\end{pf}
\end{appendix}

\begin{supplement}[id=suppA]
\stitle{Supplement to ``Batched bandit problems''\\}
\slink[doi]{10.1214/15-AOS1381SUPP} %[doi,text={...}] - jei reikia
%suskaldyti doi
\sdatatype{.pdf}
\sfilename{aos1381\_supp.pdf}
\sdescription{The supplementary material \cite{PerRigCha15supp}
contains additional simulations, including some using real data.}
\end{supplement}

%\begin{supplement}[id=suppA]
%\stitle{Supplementary Material}
%\slink[doi]{COMPLETED BY THE TYPESETTER}
%\sdescription{The supplementary material \cite{PerRigCha15supp}
%contains additional simulations, especially on real data.}
%\end{supplement}

% imsref loaded by daiva.urboniene, 2015-10-20 08:55:57

\printaddresses
\end{document}